\documentclass[preprint,12pt]{elsarticle}
\usepackage{tikz}
\usepackage{algpseudocode}
\usepackage[ruled]{algorithm}
\newtheorem{thm}{Theorem}
\newtheorem{lem}[thm]{Lemma}
\newdefinition{rmk}{Remark}
\newtheorem{example}{Example}
\newproof{pf}{Proof}
\newdefinition{definition}{Definition}

\newtheorem{corollary}{Corollary}

\usepackage{amsmath}
\usepackage{amssymb}

\DeclareMathOperator{\Pos}{Pos}

\begin{document}

\begin{frontmatter}



\title{Structural and Combinatorial Properties of 2-swap Word
Permutation Graphs}

\author[ST]{Duncan Adamson}
\author[LRC,CS]{Nathan Flaherty}
\author[CS]{Igor Potapov}
\author[CS]{Paul G. Spirakis}

\affiliation[ST]{organization={Department of Computer Science, University of St Andrews}, addressline = {North Haugh}, city={KY18 9SX}, country = {United Kingdom} }

\affiliation[LRC]{organization={Leverhulme Research Centre,University of Liverpool},%Department and Organization
            addressline={51 Oxford St}, 
            city={Liverpool},
            postcode={L7 3NY}, 
            country={United Kingdom}}

\affiliation[CS]{organization={Department of Computer Science,University of Liverpool},%Department and Organization
            addressline={Ashton Street}, 
            city={Liverpool},
            postcode={L69 3BX}, 
            country={United Kingdom}}

\begin{abstract}
In this paper, we study the graph induced by the \emph{$2$-swap} permutation (also known as a transposition) on words with a fixed Parikh vector.
%A $2$-swap is defined as a pair of positions $s = (i, j)$ where the word $w$ induced by the swap $s$ on $v$ is $v[1] v[2] \dots v[i - 1] v[j] v[i+1] \dots v[j - 1] v[i] v[j + 1] \dots v[n]$.
Informally, a $2$-swap is a permutation which swaps exactly two symbols in a word, leaving all others unchanged.   
With these permutations, we define the \emph{Configuration Graph}, $G(P)$ for a given Parikh vector.
Each vertex in $G(P)$ corresponds to a unique word with the Parikh vector $P$, with an edge between any pair of words $v$ and $w$ if there exists a 2-swap $s$ such that $v \circ s = w$.
We provide several key combinatorial properties of this graph, including the exact diameter of this graph, the clique number of the graph, and the relationships between subgraphs within this graph.
Additionally, we show that for every vertex in the graph, there exists a Hamiltonian path starting at this vertex.
Finally, we provide an algorithm enumerating these paths from a given input word of length $n$ with a delay of at most $O(\sigma \log n)$ between outputting edges, requiring $O(n \log n)$ preprocessing. 
\end{abstract}




\end{frontmatter}






\section{Introduction}

In information theory and computer science, there are several well-known edit distances between strings which are based on insertions, deletions and  substitutions of single characters or various permutations of several  characters, including swaps of adjacent or non-adjacent  characters, shuffling, etc.~\cite{Crochemore.2003,ganczorz2018edit,Maji2015}.\looseness=-1

These operations are well motivated by problems in physical science, for example, the biological swaps which occur at a gene level are non-adjacent swap operations of two symbols (mutation swap operator) representing gene mutations~\cite{AMIR2013}.
In recent work on Crystal Structure Prediction  the swap operation on a pair of symbols in a given word representing layers of atomic structures was used to generate new permutations of those layers, with the aim of exploring the configuration space of crystal structures ~\cite{mcemma}.
In computer science string-to-string correction has been studied for adjacent swaps~\cite{Levy21} and also in the context of sorting networks~\cite{ANGEL2007}, motion on graphs and diameter of permutation groups~\cite{kornhauseR1984coordinating}.
In group theory, the distance between two permutations (the Cayley distance) measures the minimum number of transpositions of elements needed to turn one into the other
\cite{KONSTANTINOVA}.\looseness=-1

% In this paper, we define the \emph{condifuration graph} as a graph mapping the set of words with a given Parikh vector to the vertices, with edges connecting vertices
A \emph{configuration graph} is a graph where words (also known as strings) are represented by vertices  and operations by edges between the strings.
For example, one may define the operations as the standard suite of edits (insertions, deletions, and substitutions), with each edge corresponding to a pair of words at an edit distance of one.
In such a graph, the distance between any pair of words corresponds to the edit distance between these words.
In this paper, we study the structural properties of such graphs defined by swap operations of two symbols on a given word (2-swap permutations), a permutation defined by a pair of indices $(i, j)$ and changing a word $w$ by substituting the symbol at position $i$ with that at position $j$, and the symbol at position $j$ with that at position $i$.
As the number of occurrences of each symbol in a given word can not be changed under this operation, we restrict our work to only those words with a given Parikh vector\footnote{
The Parikh vector of a word  $w$ denotes a vector with the number of occurrences of the symbols in the word $w$.
% The standard permutation can be seen as a permutation of the word with all distinct symbols.
}.
We focus on studying several fundamental properties of the structure of these graphs, most notably the diameter, clique number, number of cliques, and the Hamiltonicity of the graph.
Similar problems have been heavily studied for Cayley graphs~\cite{KONSTANTINOVA}, and permutation graphs~\cite{GODDARD2003}.
It has  been conjectured that the diameter of the symmetric group of degree $n$ is polynomially bounded in $n$, where only recently the exponential upper bound~\cite{babai1992diameter} was replaced by a quasipolynomial upper bound~\cite{helfgott2014diameter}. The diameter problem has additionally been studied with respect to a random pair of generators for symmetric groups~\cite{HELFGOTT2015}. In general, finding the diameter of a Cayley graph of a permutation group is NP-hard and finding the distance between two permutations in directed Cayley graphs of permutation groups is PSPACE-hard~\cite{JERRUM1985}.\looseness=-1

We also build upon previous work done into Combinatorial Gray Codes which, in general, give an ordering of some objects with two consecutive objects differ by a ``small change", The original Gray codes lists strings with subsequent strings having an edit distance of one. An extensive introduction to Combinatorial Gray codes can be found in \cite{Mtze2023}.
In \cite{aggarwal19991} Takoaka provides a combinatorial Gray code for multiset permutations (i.e. words) with constant delay however we wish to start from any word which is not a feature of their algorithm. Further, in Section 13.2.4 of \cite{Arndt2011} there is a $O(n)$ time algorithm for generating a combinatorial Gray code of multiset permutations with transpositions/2-swaps however this is provided without formal proof.
%with sections 4.1 (combinations) and 5.6 (multiset permutations) of particular relevance to our work. 
%GRAY CODES - \cite{aggarwal19991} <- constant delay but not quite the same
%\cite{Mtze2023} torsten's journal article.<- A good survey on these things can be found here, 5.6 concerns gray codes on Multiset permutations which is very similar to our construction.
%\cite{Arndt2011}<- section 13.2.4 contains an algorithm (with no proof of correctness no formal proof of correctness of this algorithm exists, and further the existing algorithm has a worst-case $O(n)$, compared with the $O(\sigma log n)$ ($\sigma$ being the number of unique elements in the set) delay we provide, making our result somewhat interesting

To develop efficient exploration strategies for these graphs it is essential to investigate its structural and combinatorial properties. As mentioned above the problem is motivated by problems arising in chemistry regarding Crystal Structure Prediction (CSP) which is computationally intractable in general~\cite{adamson2021hardness,CSP_undecidable}. In current tools~\cite{mcemma,collins2018flexible}, chemists rely on representing crystal structures as a multiset of discrete blocks, with optimisation performed via a series of permutations, corresponding to swapping blocks. Understanding reachability properties under the swap operations can help to evaluate and improve various heuristic space exploration tools and extend  related combinatorial toolbox~\cite{kCentreNecklaces,rankingBracelets}.\looseness=-1

\noindent
\subsection{Our Results}
We provide several key combinatorial properties of the  graph defined by $2$-swap permutations over a given word. First, we show that this graph is \emph{locally isomorphic}, that is, the induced subgraph of radius $r$ centred on any pair of vertices $w$ and $u$ are isomorphic. We strengthen this by providing an exact diameter on the graph for any given Parikh vector. Finally, we show that, for every vertex $v$ in the graph, there is a Hamiltonian path starting at $v$. We build upon this by providing a novel algorithm for enumerating the Hamiltonian path starting at any given vertex $v$ in a binary graph with at most $O(\log n)$ delay between outputting the swaps corresponding to the transitions made in the graph, after $O(n \log n)$ preprocessing. We extend this to general alphabets, providing an algorithm enumerating the Hamiltonian path with $O(\sigma \log n)$ delay between outputs after $O(n \log n)$ preprocessing. Our enumeration results correlate well with the existing work on the enumeration of words. This includes work on explicitly outputting each word with linear delay~\cite{adamson2023words,Ruskey1992}, or outputting an implicit representation of each word with either constant or logarithmic delay relative to the length of the words~\cite{ackerman2009three,ackerman2009efficient,makinen1997lexicographic,SchmidS21,SchmidS22}. 
The surveys~\cite{mutze2022combinatorial,wasa2016enumeration} provide a comprehensive overview of a wide range of enumeration results.

\section{Preliminaries}

Let $\mathbb{N} = \{1, 2, \dots\}$ denote the set of natural numbers, and $\mathbb{N}_0 = \mathbb{N} \cup \{ 0 \}$. We denote by $[n]$ the set $\{1, 2, \dots, n\}$ and by $[i, n]$ the set $\{i, i + 1, \dots, n\}$, for all $i, n \in \mathbb{N}_0, i \leq n$. An \emph{alphabet} $\Sigma$ is an ordered, finite set of symbols. Tacitly assume that the alphabet $\Sigma = [\sigma] = \{1, 2, \dots, \sigma\}$, where $\sigma = \vert \Sigma \vert$. We treat each symbol in $\Sigma$ both as a symbol and by the numeric value, i.e. $i \in \Sigma$ represents both the symbol $i$ and the integer $i$. A \emph{word} is a finite sequence of symbols from a given alphabet. The length of a word $w$, denoted $\vert w \vert$, is the number of symbols in the sequence. The notation $\Sigma^n$ denotes the set of $n$-length words defined over the alphabet $\Sigma$, and the notation $\Sigma^*$ denotes the set of all words defined over $\Sigma$.\looseness=-1

For $i \in [\vert w \vert]$, the notation $w[i]$ is used to denote the $i^{th}$ symbol in $w$, and for the pair $i, j \in [\vert w \vert], w[i, j]$ is used to denote the sequence $w[i] w[i + 1] \dots w[j]$, such a sequence is called a \emph{factor} of $w$. We abuse this notation by defining, for any pair $i, j \in [\vert w \vert]$ such that $j < i$, $w[i, j] = \varepsilon$, where $\varepsilon$ denotes the empty string.\looseness=-1

\begin{definition}[2-swap]
    \label{def:2-swap}
    Given a word $w \in \Sigma^{n}$ and pair $i, j \in [n], i < j$ such that $w[i] \neq w[j]$, the \emph{2-swap} of $w$ by $(i, j)$, denoted $w \circ (i, j)$, returns the word 
    \begin{center}
    $w[1, i - 1] w[j] w[i + 1, j - 1] w[i] w[j + 1, n]$.  
    \end{center}
\end{definition}

\begin{example}
    \label{ex:2-swap}
    Given the word $w = 11221122$ and pair $(2, 7)$, $w \circ (2, 7) = 12221112$.
\end{example}

Given a word $w \in \Sigma^n$, the \emph{Parikh vector} of $w$, denoted $P(w)$ is the $\sigma$-length vector such that the $i^{th}$ entry of $P(w)$ contains the number of occurrences of symbol $i$ in $w$, formally, for $i \in [\sigma]$, $P(w)[i] = \vert \{j \in [n] \mid w[j] = i \} \vert$, where $n=\vert w\vert$. For example, the word $w = 11221122$ has Parikh vector $(4,4)$. The set of words with a given Parikh vector $P$ over the alphabet $\Sigma$ is denoted $\Sigma^{*|_P}$, formally $\Sigma^{*|_P} = \{w \in \Sigma^* \mid P(w) = P \}$. Unless stated otherwise we define $n:= \sum_{i\in[\sigma]} P[i]$. It is notable that $\vert \Sigma^{*|_P} \vert = \frac{n!}{\prod_{i \in [\sigma]}P[i]!}$ since it means the configuration graph (defined below) is of exponential size in $n$.

%\subsection{Graphs}

\begin{definition}
    \label{def:configuration_grap}
    For a given alphabet $\Sigma$ and Parikh vector $P$, the \emph{configuration graph} of $\Sigma^{*|_P}$ is the undirected graph $G(P) = \{V(P), E(P) \}$ where:
    \begin{itemize}
        \item $V(P) = \{v_w \mid w \in \Sigma^{*|_P}\}$.
        \item $E(P) = \{\{v_{w}, v_u\} \in V(P) \times V(P) \mid \exists i,j \in [n]$ s.t. $w \circ (i, j) = u\}$.
    \end{itemize}
\end{definition}

\begin{figure}
    \centering
     \begin{tikzpicture}[scale=3.7]
  \tikzset{vertex/.style = {shape=circle,draw}}
      \draw
        (1, 0.0) node[vertex] (11122){11122}
        (0.809, 0.588) node[vertex] (11212){11212}
        (0.309, 0.951) node[vertex] (11221){11221}
        (-0.309, 0.951) node[vertex] (12112){12112}
        (-0.809, 0.588) node[vertex] (12121){12121}
        (-1.0, -0.0) node[vertex] (12211){12211}
        (-0.809, -0.588) node[vertex] (21112){21112}
        (-0.309, -0.951) node[vertex] (21121){21121}
        (0.309, -0.951) node[vertex] (21211){21211}
        (0.809, -0.588) node[vertex] (22111){22111};
      \begin{scope}[-]
        \draw (11122) to (21112);
        \draw (11122) to (12112);
        \draw (11122) to (11212);
        \draw (11122) to (21121);
        \draw (11122) to (12121);
        \draw[red] (11122) to (11221);
        \draw (11212) to (21112);
        \draw[red] (11212) to (12112);
        \draw (11212) to (21211);
        \draw (11212) to (12211);
        \draw[red] (11212) to (11221);
        \draw (11221) to (21121);
        \draw (11221) to (12121);
        \draw (11221) to (21211);
        \draw (11221) to (12211);
        \draw (12112) to (21112);
        \draw (12112) to (22111);
        \draw (12112) to (12211);
        \draw[red] (12112) to (12121);
        \draw (12121) to (21121);
        \draw (12121) to (22111);
        \draw[red] (12121) to (12211);
        \draw[red] (12211) to (21211);
        \draw (12211) to (22111);
        \draw[red] (21112) to (22111);
        \draw (21112) to (21211);
        \draw[red] (21112) to (21121);
        \draw (21121) to (22111);
        \draw[red] (21121) to (21211);
        \draw (21211) to (22111);
      \end{scope}
    \end{tikzpicture}
    
    \caption{The configuration graph $G(3,2)$ with Hamiltonian path shown in red.}
    \label{fig:configuration_graph}
\end{figure}
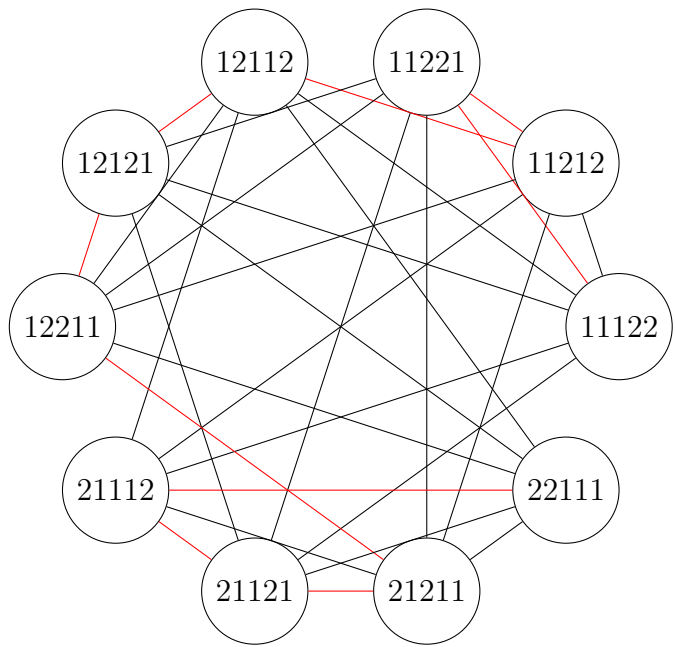

Informally, the configuration graph for a given Parikh vector $P$ is the graph with each vertex corresponding to some word in $\Sigma^{*|_P}$, and each edge connecting every pair of words $w, u \in \Sigma^{*|_P}$ such that there exists some 2-swap transforming $w$ into $u$. Figure \ref{fig:configuration_graph} provides an example of the configuration graph when $P=(3,2)$.

A \emph{path} (also called a \emph{walk}) in a graph is an ordered set of edges such that the second vertex in the $i^{th}$ edge is the first vertex in the $(i + 1)^{th}$ edge, i.e. $p = \{(v_1, v_2), (v_2, v_3), \dots, (v_{\vert p \vert}, v_{\vert p \vert} + 1)\}$. Note that a path of length $i$ visits $i + 1$ vertices. A path $p$ \emph{visits} a vertex $v$ if there exists some edge $e \in p$ such that $v \in e$. A \emph{cycle} (also called a \emph{circuit}) is a path such that the first vertex visited is the same as the last. A \emph{Hamiltonian} path $p$ is a path visiting each vertex exactly once, i.e. for every $v \in V$, there exists at most two edges $e_1, e_2 \in p$ such that $v \in e_1$ and $v \in e_2$. A cycle is Hamiltonian if it is a Hamiltonian path and a cycle. A path $p$ \emph{covers} a set of vertices $V$ if, for every $v \in V$, there exists some $e \in p$ such that $v \in e$. Note that a Hamiltonian path is a path cover of every vertex in the graph. and a Hamiltonian cycle is a cycle cover of every vertex in the graph.
We use the notation $\max_{<i} A$ to denote the largest value in $A$ which is less than $i$. And similarly $\min_{>i} A $ denotes the smallest value in $A$ which is greater than $i$.

The \emph{distance} between a pair of vertices $v, u \in V$, denoted $D(v, u)$ in the graph $G$ is the smallest value $d \in \mathbb{N}_0$ for which there exists some path $p$ of length $d$ covering both $v$ and $u$, i.e. the minimum number of edges needed to move from $v$ to $u$. If $v = u$, then $D(v,u)$ is defined as $0$. The \emph{diameter} of a graph $G$ is the maximum distance between any pair of vertices in the graph, i.e. $\max_{v, u \in V} D(v, u)$.\looseness=-1

Given two graphs $G = (V, E)$ and $G' = (V', E')$, $G$ is \emph{isomorphic} to $G'$ if there exists a bijective mapping $f : V \mapsto V'$ such that, for every $v, u \in V$,  $(v, u) \in E$ if and only if $(f(v), f(u)) \in E'$. The notation $G \cong G'$ is used to denote that $G$ is isomorphic to $G'$, and $G \not\cong G'$ to denote that $G$ is not isomorphic to $G'$. A \emph{subgraph} of a graph $G = (V, E)$ is a graph $G' = (V', E')$ such that $V' \subseteq V$ and $E' \subseteq E$. A \emph{clique} $G'=(V',E')$ is a subgraph, $G' \subseteq G$ which is complete (i.e. for all $u,v \in G',$  $(u,v) \in E'$). And the clique number $\omega$ of a graph $G$ is the size of the largest clique in $G$.\looseness=-1

\section{Basic Properties of the Configuration Graph}
\label{sec:configuration_graph}

In this section, we provide a set of combinatorial results on the configuration graph. We first show that every subgraph, $G_r(v)$, of the configuration graph $G(P) = (V, E)$ with the vertex set $V'(v) = \{u \in V \mid D(v, u) \leq \ell\}$ for $v \in V$, and edge set $E' = (V' \times V') \cap E$ are isomorphic. We build on this by providing a tight bound on the diameter of these graphs.
We start by considering some local structures within the graph.\looseness=-1

\begin{lem}\label{lem:cliques}
Given a Parikh vector $P$ with associated configuration graph $G(P)$, each vertex $v \in V(P)$ belongs to $\sum_{j \in \Sigma}\prod_{i \in \Sigma \setminus \{j\}} P[i]$ maximal cliques, with the size of each such clique being in $\{P[i]+1 \mid i \in \Sigma\}$.\looseness=-1
\end{lem}

\begin{pf}
Consider first the words with Parikh vector $P = (k, 1)$. Note that every word in $\Sigma^{*|_P}$ consists of $k$ copies of the symbol $1$, and one copy of the symbol $2$. Therefore, given any pair of words $w, u \in \Sigma^{*\vert_P}$ such that $w[i] = u[j] = 2$, the 2-swap $(i, j)$ transforms $w$ into $u$ and hence there exists some edge between $w$ and $u$.
Hence $G(P)$ must be a complete graph of size $k + 1$.\looseness=-1

In the general case, consider the word $w\in \Sigma^{*|_P}$ where $P = (k_1$, $k_2$, $\dots$, $k_{\sigma})$, and $n = \sum_{k_i \in P} k_i$. Let $\Pos(w, i) = \{j \in [n] \mid w[j] = i\}$. Given $i,j \in [\sigma], i \neq j$, let $i_1, i_2 \in \Pos(w, i)$ and $j_1 \in \Pos(w, j)$ be a set of indices. Let $v_1 = w \circ (i_1, j_1)$ and $v_2 = w \circ (i_2, j_1)$. Then, $v_1[i_1] = v_2[i_2], v_2[i_2] = v_1[i_1]$, and $v_1[j_1] = v_2[j_1]$. Further, for every $\ell \in [n]$ such that $\ell \notin \{i_1, i_2, j_1\}$, $v_1[\ell] = v_2[\ell]$ as these positions are unchanged by the swaps. Therefore, $v_1 = v_2 \circ (i_1, i_2)$, and hence these words are connected in $G(P)$. Further, as this holds for any $j_1 \in \Pos(v, j)$, the set of words induced by the swaps $(j,\ell)$, for some fixed $\ell \in \Pos(w, i)$ correspond to a clique of size $P[j] + 1$. Therefore, there exists $\prod_{i \in \Sigma \setminus \{j\}} P[i]$ cliques of size $P[j] + 1$ including $w$, for any $j \in \Sigma$.\looseness=-1

We now show that the cliques induced by the set of swaps $S(i, j) = \{(i', j) \mid i' \in \Pos(w, w[i])\}$ are maximal. Let $C(i, j, w) = \{ w \} \cup \{w \circ (i', j) \mid (i', j) \in S(i, j) \}$, i.e. the clique induced by the set of swaps in $S(i, j)$. Consider a set of swaps, $(i_1, j_1), (i_2, j_1), (i_1, j_2)$ and $(i_2, j_2)$, where $i_1, i_2 \in \Pos(w, i)$ and $j_1, j_2 \in \Pos(w, j)$. Let $v_{1, 1} = w \circ (i_1, j_1), v_{2, 1} = w \circ (i_2, j_1), v_{1, 2} = w \circ (i_1, j_2)$ and $v_{2,2} = w \circ (i_2, j_2)$. Note that $\{w, v_{1, 1}, v_{2, 1}\} \subseteq C(i, j_1), \{w, v_{1, 1}, v_{1, 2}\} \subseteq C(j, i_1)$, $\{w, v_{2, 1}, v_{2, 2}\} \subseteq C(i, j_2)$ and $\{w, v_{2, 1}, v_{2, 2}\} \subseteq C(j, i_2)$.\looseness=-1

We now claim that there exists no swap transforming $v_{1, 1}$ in to $v_{2, 2}$. Observe first that, for every $\ell \in [\vert w \vert]$ such that $\ell \notin \{i_1, i_2, j_1, j_2\}$, $v_{1, 1}[\ell] = v_{2, 2}[\ell]$. As $v_{1, 1}[i_1] = v_{2, 2}[j_1]$, and $v_{1, 1}[j_1] = v_{2, 2}[i_2]$, exactly two swaps are needed to transform $v_{1, 1}$ into $v_{2, 2}$. Therefore, for any pair of swaps $(i_1, j_1), (i_2, j_2) \in \Pos(i, w) \times \Pos(j, w)$, such that $i_1 \neq i_2$ and $j_1 \neq j_2$, the words $w \circ (i_1, j_1)$ and $w \circ (i_2, j_2)$ are not adjacent in $G(v)$. Similarly, given a set of indices $i' \in \Pos(i, w), j' \in \Pos(j, w)$ and $\ell' \in \Pos(\ell, w)$ and swaps $(i', j'), (i', j')$, observe that as $w[j'] \neq w[\ell']$, the distance between $w \circ (i', j')$ and $w \circ (i', \ell')$ is $2$. Therefore, every clique induced by the set of swaps $S(i, j) = \{(i', j) \mid i' \in \Pos(w, w[i])\}$ is maximal.\looseness=-1
\end{pf}

\begin{corollary}
\label{col:only_belongs_to_the_cliques}
    Let $v \in V(P)$ be the vertex in $G(P)$ corresponding to the word $w \in \Sigma^{*|_P}$. Then, $v$ belongs only to the maximal cliques corresponding to the set of words $\{w \circ (i, j) \mid i \in \Pos(w, x)\}$ for some fixed symbol $x \in \Sigma$ and position $j \in [\vert w \vert], w[j] \neq x$, where $\Pos(w, x) = \{i \in [\vert w \vert] \mid, w[i] = x\}$.\looseness=-1
\end{corollary}

Now since the edges in each maximal clique only swap two types of symbols we have the following corollary for the number of cliques.

\begin{corollary}
    \label{col:number_of_cliques}
    There are $\sum_{(i,j) \in \Sigma \times \Sigma }  \frac{\left( \sum_{k \in \Sigma \setminus \{i,j\}} P[k]\right )!}{\prod_{k\in \Sigma \setminus \{i,j\}} P[k]!} $ maximal cliques in $G(p)$. 
\end{corollary}

\begin{corollary}
The clique number $\omega(G(P))$ is equal to $\max_{i\in \Sigma} P[i]  +1$.\looseness=-1
\end{corollary}

\begin{lem}\label{lem:isomorphic_subgraphs}
    Let $G_r(v)$ be the subgraph of $G(P)$ induced by all vertices of distance at most $r$ away from a given vertex $v$. Then, for any pair of vertices $u,v \in V$ and given any $r \in \mathbb{Z}^+$, $G_r(u) \cong G_r(v)$.\looseness=-1
\end{lem}

\begin{pf}
Let $\pi \in S_n$ be the permutation such that $u \circ \pi = v$. We use the permutation $\pi$ to define an isomorphism $f: G_r(u) \rightarrow G_r(v)$ such that $f(w) = w \circ \pi $. In order to show that $f$ is an isomorphism we need to show that it preserves adjacency. We start by showing that for every word, $w \in G_1(u),  f(w) \in G_1(v)$.\looseness=-1

Let $\tau = (\tau_1,\tau_2)$ be the 2-swap such that $w = u \circ \tau$. We now have 3 cases for how $\pi$ and $\tau$ interact, either none of the indices in $\tau$ are changed by $\pi$, just one of $\tau_1$ or $\tau_2$ are changed by $\pi$, or both $\tau_1$ and $\tau_2$ are changed by $\pi$. In the first case, $f(w)$ is adjacent to $v$ as $v \circ \tau = f(w)$. In the second case, let $(\tau_1, \tau_2)$ be a swap that that $\pi[\tau_1] = \tau_1$, i.e. $\tau_1$ is not changed by the permutation $\pi$. We define a new swap $\tau'$ such that $v \circ \tau' = f(w)$. Let $x, y \in [n]$ be the positions in $v$ such that $\pi[x] = \tau_2$ and $\pi[\tau_2] = y$. Now, let $\tau' = (\tau_1, y)$. Observe that $v[y] = w[\tau_2]$, and $v[\tau_1] = w[\tau_1]$. Therefore, the word $v \circ \tau' = u \circ \tau \circ \pi$. Note that as the ordering of the indices in the swap does not change the swap, the same argument holds for the case when $\pi[\tau_2] = \tau_2$. In the final case, let $\tau' = (\pi[\tau_1], \pi[\tau_2])$. Note that by arguments above, $u[\pi[\tau_1]] = v[\tau_1]$ and $u[\pi[\tau_2]] = v[\tau_2]$,  and hence $v \circ \tau' = u \circ \tau \circ \pi$. Repeating this argument for each word at distance $\ell \in [1, r]$ proves this statement.\looseness=-1
\end{pf}

We now provide the exact value of the diameter of any configuration graph $G(P)$. Theorem \ref{thm:diameter} states the explicit diameter of the graph, with the remainder of the section dedicated to proving this result.\looseness=-1

\begin{thm}\label{thm:diameter}
    The diameter of the Configuration Graph, $G(P)$ for a given Parikh vector $P$ is $n-\max_{i \in \Sigma} P[i]$.\looseness=-1
\end{thm}

Theorem \ref{thm:diameter} is proven by first showing that the upper bound matches $n-\max_{i \in \Sigma} P[i]$ (Lemma \ref{lem:diameter_upper_bound}). We then show that the lower bound on the diameter matches the upper bound (Lemma \ref{lem:diameter_lower_bound}), concluding our proof of Theorem \ref{thm:diameter}.\looseness=-1

\begin{lem}[Upper Bound of Diameter]
    \label{lem:diameter_upper_bound}
      The diameter of the Configuration Graph, $G(P)$ for a given Parikh vector $P$ %, with 2-swaps as edges
      is at most $n-\max_{i \in \Sigma} P[i]$.\looseness=-1
\end{lem}

\begin{pf}[Proof of Upper Bound]
    This claim is proven by providing a procedure to determine a sequence of $n-\max_{i \in [\sigma]} P[i]$  swaps to transform any word $w \in \Sigma^{*|_P}$ into some word $v \in \Sigma^{*|_P}$. We assume, without loss of generality, that $P[1] \geq P[2] \geq \dots \geq P[\sigma]$. \looseness=-1
    The procedure described by Algorithm \ref{alg:paths} operates by iterating over the set of symbols in $\Sigma$, and the set of occurrences of each symbol in the word. At each step, we have a symbol $x \in [2, \sigma]$ and index $k \in [1, P[x]]$. The procedure finds the position $i$ of the $k^{th}$ appearance of symbol $x$ in $w$, and the position $j$ of the $k^{th}$ appearance of $x$ in $v$. Formally, $i$ is the value such that $w[i] = x$ and $\vert \{i' \in [1, i - 1] \mid w[i'] = x\} \vert = k$ and $j$ the value such that $v[j] = x$ and $\vert \{j' \in [1, j - 1] \mid v[j'] = x\} \vert = k$. Finally, the algorithm adds the swap $(i, j)$ to the set of swaps, and then moves to the next symbol.\looseness=-1
    
   \begin{algorithm}
    \caption{Procedure to select 2-swaps to generate a path from $w$ to $v$.}
    \begin{algorithmic}[1]
        \State{$S \leftarrow \emptyset$} \Comment{Set of 2-swaps}
        \For{$x \in \Sigma \setminus \{1\} $}
        \For{$1 \leq k \leq P_x$}
        \State $i \leftarrow$ index of $k^\textrm{th}$ occurrence of $x$ in $w$
        \State $j \leftarrow$ index of $k^\textrm{th}$ occurrence of $x$ in $v$
        \State $S \leftarrow S \cup (i,j)$
        \EndFor
        \State{Apply all 2-swaps in $S$ to $w$ and set $S \leftarrow \emptyset$} 
        \EndFor
    \end{algorithmic}
    \label{alg:paths}
\end{algorithm}

   This procedure requires one swap for each symbol in $w$ other than $1$, giving a total of $n - \max_{i \in [\sigma]} P[i]$ swaps. Note that after each swap, the symbol at position $j$ of the word is the symbol $v[j]$. Therefore, after all swaps have been applied, the symbol at position $j \in \{i \in [1, \vert w \vert] \mid v[i] \in \Sigma \setminus \{1 \}\}$ must equal $v[j]$. By extension, for any index $i$ such that $v[i] = 1$, the symbol at position $i$ must be $1$, and thus equal $v[i]$. Therefore this procedure transforms $w$ into $v$.
\end{pf}

In order to prove the lower bound on the diameter (i.e. that $\texttt{diam}(G) \geq n- \max_{i \in [\sigma]} P[i]$) we introduce a new auxiliary structure, the \emph{2-swap graph}. Informally, the 2-swap graph, defined for a pair of words $w, v \in \Sigma^{*|_P}$ and denoted $G(w, v) = \{V(w, v), E(w, v)\}$ is a directed graph such that the edge $(u_i, u_j) \in V(w, v) \times V(w, v)$ if and only if $w[i] = v[j]$. Note that this definition allows for self-loops.\looseness=-1

\begin{definition}[2-swap Graph]
Let $w, v \in \Sigma^{*|_P}$ be a pair of words. The \emph{2-swap graph} $G(w, v) = \{V(w, v), E(w, v)\}$ contains the vertex set $V(w, v) = \{u_1, u_2, \dots, u_{\vert  w \vert}\}$ and edge set $E(v,w) = \{(v_i, v_j) \in V(w, v) \times V(w, v) \mid w[i] = v[j]$ or $v[i] = w[j]\}$. The edge set, $E$, is defined as follows, for all $i,j \in \Sigma$ there exists an edge $(i,j) \in E$ if and only if $w[i] = v[j]$.\looseness=-1
\end{definition}

\begin{figure}[t]
\centering
\begin{tikzpicture}[scale=1.3]
    \tikzset{vertex/.style = {shape=circle,draw,minimum size=2em}}
    \tikzset{edge/.style = {->}}
    
    \node[vertex](1) at (0,0){1};
    \node[vertex](2) at (2,0){2};
    \node[vertex](3) at (3.9,0){3};
    \node[vertex](4) at (6.1,0){4};
    \node[vertex](5) at (8,0){5};
    \node[vertex](6) at (10,0){6};
    
    \draw[edge]  (1) to [bend left=38] node[above] {a} (4); 
    \draw[edge](4) to [bend left=44] node[above] {b} (1); 
    \draw[edge] (1) to  [bend left=40]node[above] {a} (5); 
    \draw[edge] (1) to [bend left=46] node[above] {a}  (6); 
    \draw[edge] (2) to [bend left=28] node[above] {a} (4); 
    \draw[edge] (2) to [bend left=35]node[above] {a}  (5); 
    \draw[edge] (5) to [bend left=48] node[above] {b} (1); 
    \draw[edge] (2) to [bend left=40] node[above] {a} (6); 
    \draw[edge] (6) to [bend left=35] node[above] {c} (2); 
    \draw[edge] (3) to  node[above] {a} (4); 
    \draw[edge] (4) to [bend left=15] node[below] {b} (3); 
    \draw[edge] (3) to [bend left=28] node[above] {a} (5); 
    \draw[edge] (5) to [bend left=38] node[above] {b} (3); 
    \draw[edge] (3) to [bend left=37] node[above] {a} (6); 
    
\end{tikzpicture}
\caption{The graph $G(aaabbc,bcbaaa)$ with the edge $(i,j)$ labelled by the symbol $w[i] (= v[j])$.}

\label{fig:graphwv}
\end{figure}
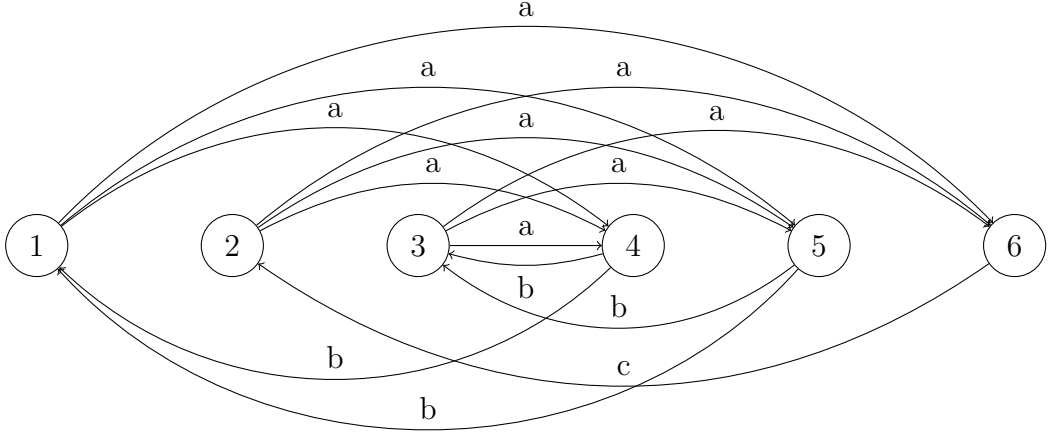

An example of the 2-swap graph is given in Figure \ref{fig:graphwv}.

\begin{lem}
Let $G(w,v)$ be a graph constructed as above for transforming $w$ into $v$ using 2-swaps.
Then, there exists a procedure to convert any cycle cover of $G(w,v)$, $\mathcal{C}$, into a $w-v$ path in $G(P)$\looseness=-1
%Then, any cycle cover $\mathcal{C}$ of $G(w,v)$ can be converted into a set of 2-swaps to transform $w$ into $v$.\looseness=-1
\end{lem}

\begin{pf}
Let $C \in \mathcal{C}$ be a cycle where $C = (e^1, e^2, \dots, e^{|C|})$ and $e^i_2 = e^{i + 1 \bmod |C|}_1$. The $w-v$ path (i.e. a sequence of 2-swaps) is constructed as follows. Starting with $i = 1$ in increasing value of $i \in [|C| - 1]$, the 2-swap $(e^1_1, e^i_2)$ is added to the set of 2-swaps $S$. Where $e^i_1$ and $e^i_2$ are the endpoints of edge $e^i$ for each $i$. \looseness=-1

Assume, for the sake of contradiction, that $S$ does not correspond to a proper set of 2-swaps converting $w$ into $v$. Then, there must exist some symbol at position $i$ such that the symbol $w[i]$ is placed at some position $j$ such that $w[i] \neq v[j]$. As $w_i$ must be placed at some position that is connected to node $i$ by an edge, there must be an edge between $i$ and $j$, hence $w_i = v_j$, contradicting the construction of $G(w, v)$. Therefore, $S$ must correspond to a proper set of 2-swaps.\looseness=-1
\end{pf}

\begin{corollary}
\label{col:cycles_to_total_2-swaps}
Let $\mathcal{C}$ be a cycle cover of $G(w,v)$. Then there exists a set of $\sum_{c \in \mathcal{C}} |c| - 1$ 2-swaps transforming $w$ in to $v$.\looseness=-1
\end{corollary}

\begin{corollary}
\label{col:disjoint}
Let $S$ be the smallest set of 2-swaps transforming $w$ in to $v$, then $S$ must correspond to a vertex disjoint cycle cover of $G(w,v)$.\looseness=-1
\end{corollary}

\begin{pf}
For the sake of contradiction, let $S$ be the smallest set of 2-swaps transforming $w$ into $v$, corresponding to the cycle cover $\mathcal{C}$ where $\mathcal{C}$ is not vertex disjoint. Let $c_1, c_2 \in \mathcal{C}$ be a pair of cycles sharing some vertex $u$. Then, following the construction above, the symbol $w_u$ must be used in two separate positions in $v$, contradicting the assumption that $v$ can be constructed from $w$ using 2-swaps. Hence $S$ must correspond to a vertex disjoint cycle cover.\looseness=-1
\end{pf}

\begin{corollary}
\label{col:cycle_cover_condition}
Given a pair of words $w,v \in \Sigma^{*|_P}$, the minimum set of 2-swaps transforming $w$ into $v$ $S$ corresponds to the vertex disjoint cycle cover of $G(w,v)$ maximising the number of cycles.\looseness=-1
\end{corollary}

\begin{lem}[Lower Bound]
\label{lem:diameter_lower_bound}
  The diameter of $G(p)$ is at least $n-\max_{i\in\Sigma}P[i]$.
\end{lem}
\begin{pf}
We assume 
w.l.o.g. that $P[1] \geq P[2] \geq \dots \geq P[\sigma]$.
Let $w,v \in \Sigma^{*|_P}$ satisfy:\looseness=-1
\vspace{-0.3cm}
$$w = (1 2 3 \dots \sigma)^{P[\sigma]} (1 2 3 \dots \sigma - 1)^{P[{\sigma - 1}]- P[\sigma]} \dots 1^{P[1] - P[2]}$$
and
\[ =  (2 3 \dots \sigma 1)^{P[\sigma]} (2 3 \dots (\sigma - 1) 1)^{P[{\sigma - 1}] - P[\sigma]} \dots 1^{P[1] - P[2]},\]

i.e. $w$ is made up of $P[1]$ subwords, each of which are of the form $12...k$, and $v$  is made up of the same subwords as $w$ but each of them has been cyclically shifted by one (for example when $P=(3,2,1)$ we have $w=123121$ and $v=231211$).
Following Corollary \ref{col:cycle_cover_condition}, the minimum number of 2-swaps needed to convert $w$ into $v$ can be derived from a vertex disjoint cycle cover of $G(w,v)$ with the maximum number of cycles.

Observe that any occurrence of symbol $\sigma$ must have an outgoing edge in $G(w,v)$ to symbol $1$, and an incoming edge from symbol $\sigma - 1$. Repeating this logic, each instance of $\sigma$ must be contained within a cycle of length $\sigma$. Removing each such cycles and repeating this argument gives a set of $P[1]$ cycles, with $P[\sigma]$ cycles of length $\sigma$, $P[{\sigma- 1}] - P[\sigma]$ cycles of length $\sigma - 1$, and generally $P[i] - P[{i + 1}]$ cycles of length $i$. This gives the number of 2-swaps needed to transform $w$ to $v$ being a minimum of $n-P[1]$ = $n-\max_{i\in\Sigma}P[i]$.
\end{pf}

Theorem \ref{thm:diameter} follows from Lemmas \ref{lem:diameter_upper_bound} and \ref{lem:diameter_lower_bound}.

\section{Hamiltonicity}

In this section, we prove that the configuration graph contains Hamiltonian paths and we provide an efficient algorithm for enumerating the vertices of a Hamiltonian path. We first show that every configuration graph of a Parikh Vector over a binary alphabet is Hamiltonian. This is then generalised to alphabets of size $\sigma$, using the binary case to build Hamiltonian paths with alphabets of size $\sigma$.\looseness=-1

\noindent
\subsection{Binary Alphabets}
For notational conciseness, given a symbol $a$, in a binary alphabet $\Sigma$, the notation $\overline{a}$ is used to denote $\overline{a} \in \Sigma, a \neq \overline{a}$, i.e. if $a = 1$, then $\overline{a} = 2$. We prove Hamiltonicity via a recursive approach that forms the basis for our enumeration algorithm. Our proof works by taking an arbitrary word in the graph $w$, and constructing a path starting with $w$. At each step of the path, the idea is to find the shortest suffix of $w$ such that both symbols in $\Sigma$ appear in the suffix. Letting $w = p s$, the path is constructed by first forming a path containing every word $p s'$, for every $s' \in \Sigma^{* \vert P(s)}$, i.e. a path from $w$ transitioning through every word formed by maintaining the prefix $p$ and permuting the suffix $s$. Once every such word has been added to the path, the algorithm repeats this process by performing some swap of the form $(\vert p \vert, i)$ where $i \in [\vert p \vert + 1, \vert w \vert]$, i.e. a swap taking the last symbol in the prefix $p$, and replacing it with the symbol $\overline{w[\vert p \vert]}$ from some position in the suffix.

This process is repeated, considering increasingly long suffixes, until every word has been covered by the path. Using this approach, we ensure that every word with the same prefix is added to the path first, before shortening the prefix. The algorithm \textsc{HamiltonEnumeration} outlines this logic within the context of the enumeration problem, where each transition is output while constructing the path.

\begin{thm}
    \label{thm:binary_graph_is_hamiltonian}
    For every Parikh vector $P \in \mathbb{N}_0^2$ and word $w \in \Sigma^{*|_P}$, there exists a Hamiltonian path starting at $w$ in the configuration graph $G(P) = (V(P), E(P))$.
\end{thm}

As stated in~\cite{buck1984gray} the binary reflected Gray code gives an ordering for the words over a Binary alphabet restricted to a certain Parikh Vector ($k,n-k$) with a single 2-swap between each subsequent word. This does indeed prove Theorem \ref{thm:binary_graph_is_hamiltonian} by providing a Hamiltonian Circuit for a given binary Parikh vector (it is worth noting that following a Hamiltonian Circuit starting at $w$ gives a Hamiltonian path from $w$.) However, we also present our own inductive proof for this case to provide an explanation of our enumeration algorithm.

\begin{pf}
    We prove this statement in a recursive manner. As a base case, consider the three vectors of length $2$ as our Parikh Vector, namely $(2, 0), (0, 2)$ and $(1, 1)$. Note that there exists only a single word with the Parikh vectors $(2, 0)$ or $(0, 2)$, and thus the graph must, trivially, be Hamiltonian. For the Parikh vector $(1,1)$, there exists only the words $1 2$ and $2 1$, connected by the 2-swap $(1,2)$ and therefore is also a Hamiltonian path and it can be found starting at either word.

    In the general case, assume that for every Parikh vector $P' = (P_1', P_2')$ with $P_1' + P_2' < \ell$, the graph $G(P')$ contains a Hamiltonian path, and further there exists such a path starting at every word in $\Sigma^{* \vert_{P'}}$. Now, let $P = (P_1, P_2)$ be an arbitrary Parikh vector such that $P_1 + P_2 = \ell$. Given some word $w \in \Sigma^{*|_P}$, observe that there must exist some Hamiltonian path starting at the word $w[2, \ell]$ in the subgraph $G'(P) = (V'(P), E'(P))$ where $V'(P) = \{u \in V(P) \mid u[1] = w[1]\}$ and $E'(P) = (V'(P) \times V'(P)) \cap E(P)$. Let $w'$ be the last word visited by the Hamiltonian path in $G'(P)$, and let $i$ be some position in $w'$ such that $w'[i] = \overline{w[1]}$. Note that there must exist Hamiltonian path starting at $(w \circ (1, i))[2, \ell]$ in the subgraph $G''(P) = (V''(P), E''(P))$ where $V''(P) = \{u \in V(P) \mid u[1] = \overline{w[1]}\}$ and $E''(P) = (V''(P) \times V''(P)) \cap E(P)$. As every vertex in $G(P)$ is either in the subgraph $G'(P)$ or $G''(P)$, the Hamiltonian paths starting at $w$ in $G'(P)$ and at $w' \circ (1, i)$ in $G''(P)$ cover the complete graph. Further, as these paths are connected, there exists a Hamiltonian path starting at the arbitrary word $w \in \Sigma^{*|_P}$, and therefore the Theorem holds.
\end{pf}

\noindent
{\bf Enumeration.}
We now provide our enumeration algorithm. Rather than output each word completely, we instead maintain the current state of the word in memory and output the swaps taken at each step, corresponding to the edges traversed in the path. This way, at any given step the algorithm may be paused and the current word fully output, while the full path can be reconstructed from only the output.
There are two key challenges behind this algorithm. First is the problem of deciding the next swap to be taken to move from the current word in the graph to the next word. Second, is the problem of minimising the worst-case delay in the output of these swaps, keeping in mind that the output is of constant size.

\noindent
\begin{algorithm}
    \begin{algorithmic}[1]
        \State \textbf{Global Variables:}
        \State Word $w \in \Sigma^n$
        \State Balanced Binary Search Tree $T_1$
        \State Balanced Binary Search Tree $T_2$
        \Function{HamiltonianEnumeration}{$(P_1, P_2) \in \mathbb{N}_0 \times \mathbb{N}_0$, pointer $last\_state$}
        \If{$P_1 = 0$ or $P_2 = 0$}
            {\textbf{Return To} $last\_state$}
        \ElsIf{$(P_1, P_2) = (1, 1)$}
            \State \textbf{Output:} $(n - 1, n)$
            \State \textbf{Remove$(T_{w[n - 1]}, n - 1)$}, \textbf{Insert$(T_{w[n - 1]}, n)$}
            \State \textbf{Remove$(T_{w[n]}, n)$}, \textbf{Insert$(T_{w[n]}, n - 1)$}
            \State $w \gets w \circ (n - 1, n)$
            \State \textbf{ReturnTo} $last\_state$
        \Else $\;$ \% Note that $P_1 + P_2 \geq 3$
                \State $R_1 \gets \max(T_1)$
                \State $R_2 \gets \max(T_2)$
                \For{$i \in \min(R_1, R_2), \min(R_1, R_2) - 1, \dots, n - (P_1 + P_2 - 1)$}
                    % \State $j \gets \max_{j' \in [i + 1, n]} T_{\overline{w[i]}}$
                    \State $j \gets \min_{j' \in [i + 1, n]} T_{\overline{w[j']}}$
                    
                    \State \textbf{Output:} $(i, j)$
                    \State \textbf{Remove}$(T_{w[i]}, i)$, \textbf{Insert}$(T_{w[i]}, j)$
                    \State \textbf{Remove}$(T_{w[j]}, j)$, \textbf{Insert}$(T_{w[j]}, i)$
                    \State $w \gets w \circ (i, j)$
                    \State $P' \gets (P_1, P_2) - P(w[i])$
                    \State HamiltonianEnumeration$(P', CurrentState())$
                \EndFor
            \State \% As the Parikh vector must have at least one value for each symbol, there is a valid swap from $n - (P_1 + P_2 - 1)$ to $i$, for some $i > n - (P_1 + P_2 - 1)$
            %\State $j \gets \min_{i \in [n - (P_1 + P_2 - 1), n]} T_{\overline{w[n - (P_1 + P_2 - 1)]}}$
            \State $j \gets \min_{i \in [n - (P_1 + P_2 - 1), n]} T_{\overline{w[i]}}$
            \State \textbf{Output} $(n - (P_1 + P_2 - 1), j)$
            \State \textbf{Remove}$(T_{w[n - (P_1 + P_2 - 1)]}, n - (P_1 + P_2 - 1))$, \textbf{Insert}$(T_{w[n - (P_1 + P_2 - 1)]}, j)$
            \State \textbf{Remove}$(T_{w[j]}, j)$, \textbf{Insert}$(T_{w[j]}, n - (P_1 + P_2 - 1))$
            \State $w \gets w \circ (n - (P_1 + P_2 - 1), j)$
            \State $P' \gets (P_1, P_2) - P(w[n - (P_1 + P_2 - 1)])$
            \State HamiltonianEnumeration$(P', last\_state)$ \% Note that this skips over the current state when returning
        \EndIf
        \EndFunction
    \end{algorithmic}
    \caption{Algorithm for enumerating a Hamiltonian path in the configuration graph defined by a Parikh vector $P$. Note that the pointer $last\_state$ is used to return at the state in the stack at position $last\_state$ to avoid needless recursion.
    The function $CurrentState()$ is used to get the current state in the stack.
    The trees $T_1$ and $T_2$ are balanced binary search trees such that every node in $T_1$ corresponds to a position of symbol $1$ in the current state of $w$, and every node in $T_2$ corresponds to a position of the symbol $2$ in $w$.
    }
    \label{alg:hamiltonian_enumeration}
\end{algorithm}

\emph{High-Level Idea.}

\noindent
From a given word $w$ with Parikh vector,$P$, the algorithm works by first finding the shortest suffix $s$ of $w$ such that there exists some pair of indices $i, j$ for which $s[i] \neq s[j]$. Using this suffix and letting $w = u s$, we find a path through every vertex in $G(P)$ with the prefix $u$. Note that following the same arguments as Theorem \ref{thm:binary_graph_is_hamiltonian}, such a path must exist. Once every word in $G(P)$ with the prefix $u$ has been visited by the path, the algorithm then enumerates every word with the prefix $u[1, \vert u \vert - 1]$, extending the current path. When adding every word with the prefix $u[1, \vert u \vert - 1]$ to the path, note that every word with the prefix $u$ has already been added, thus all that is left is to add those words with the prefix $u[1, \vert u \vert - 1] \overline{u[\vert u \vert]}$, which is achieved via the same process as before.\looseness=-1

The swaps are determined as follows. From the initial word $w$, let $R_1$ be the last occurrence of the symbol $1$ in $w$, and let $R_2$ be the last occurrence of $2$ in $w$. The first swap is made between $\min(R_1, R_2)$ and $\min(R_1, R_2) + 1$, with the algorithm then iterating through every word with the Parikh vector $P\left[w[\min(R_1, R_2), \vert w \vert]\right] - P\left[\overline{w[\min(R_1, R_2)]}\right]$.\looseness=-1

In the general case, a call is made to the algorithm with a Parikh vector $P=(P_1, P_2)$, with the current word $w$ fixed, and the assumption that no word with the prefix $w[1, \vert w \vert - (P_1 + P_2)]$ has been added to the path other than $w$. The algorithm, therefore, is tasked with iterating through every word with the current prefix. Let $R_1$ be the last occurrence of the symbol $1$, and $R_2$ be the last occurrence of the symbol $2$ in the current word. The algorithm first enumerates every word with the prefix $w[1, \min(R_1, R_2) - 1]$. Noting that there exists only a single word with the prefix $w[1, \min(R_1, R_2)]$, it is sufficient to only enumerate through those words with the prefix $w[1, \min(R_1, R_2) - 1]$ $\overline{w[\min(R_1, R_2)]}$. The first swap made by this algorithm is \newline$(\min(R_1, R_2), \min(R_1, R_2) + 1)$, allowing a single recursive call to be made to \newline$HamiltonianEnumeration(P(w \circ (\min(R_1, R_2), \min(R_1, R_2) + 1))[\min(R_1, R_2) + 1, \vert w \vert]$, \textsc{current call}$)$, where \textsc{current call} denotes the pointer to the current call on the stack. 
%Note that 
From this call the algorithm enumerates every word with the prefix $w[1, \min(R_1, R_2) - 1]$ $\overline{w[\min(R_1, R_2)]}$. As every word with the prefix  $w[1,\min(R_1, R_2)]$ has already been output and added to the path, once this recursive call has been made, every word with the prefix $w[1, \min(R_1, R_2) - 1]$ will have been added to the path. Note that the word $w$ is updated at each step, ending at the word $w'$.

After every word with the prefix $w'[1, \min(R_1, R_2) - 1]$ has been added to the path, the next step is to add every word with the prefix $w'[1, \min(R_1, R_2) - 2]]$ to the path. As every word with the prefix $w[1, \min(R_1, R_2) - 1]$ is already in the path, it is sufficient to add just those words with the prefix $w'[1, \min(R_1, R_2) - 2] \overline{w[\min(R_1, R_2)]}$ to the path. This is achieved by making the swap between $\min(R_1, R_2) - 2$, and the smallest value $i > \min(R_1, R_2) - 2$ such that $w[i] \neq  w'[\min(R_1, R_2) - 1]$, then recursively enumerating every word with the prefix $w'[1, \min(R_1, R_2) - 2]$. This process is repeated in decreasing prefix length until every word has been enumerated.

To efficiently determine the last position in the current word $w$ containing the symbols $1$ and $2$, a pair of balanced binary search trees are maintained. The tree $T_1$ corresponds to the positions of the symbol $1$ in $w$, with each node in $T_1$ being labelled with an index and the tree sorted by the value of the labels. Analogously, tree $T_2$ corresponds to the positions of the symbol $2$ in $w$. Using these trees, note that the last position in $w$ at which either symbol appears can be determined in $O(\log n)$ time, and further each tree can be updated in $O(\log n)$ time after each swap.

\begin{lem}
    \label{lem:alg_he_visits_every_node}
    Let $P$ be a Parikh vector of length $n$, and let $w \in \Sigma^{*|_P}$ be a word. % Algorithm \ref{alg:hamiltonian_enumeration} 
    \textsc{HamiltonianEnumeration} outputs a path visiting every word in $\Sigma^{*|_P}$ starting at $w$.
\end{lem}
\begin{pf}
    This lemma is proven via the same tools as Theorem \ref{thm:binary_graph_is_hamiltonian}. Explicitly, we show first that the algorithm explores every suffix in increasing length, relying on the exploration of suffixes of length $2$ as a base case, then provide an inductive proof of the remaining cases. We assume that the starting word has been fully output as part of the precomputation. With this in mind, note that there are two cases for length 2 prefixes, either the suffix contains two copies of the same symbol or one copy of each symbol. In the first case, as $w$ has been output, so has every permutation of the length 2 prefix of $w$. Otherwise, the algorithm outputs the swap $(n - 1, n)$ and returns to the previous call.

    In the general case, we assume that for some $\ell \in [n]$, every permutation of $w[n - \ell + 1, n]$ has been visited by the path. Further, we assume the algorithm can, given any word $v$, visit every word of the form $v[1, n - \ell] u$, for every $u \in \Sigma^{* | _{P(v[n - \ell + 1, n])}}$, i.e. the algorithm is capable of taking any word $v$ as an input, and visiting every word with the same Parikh vector $P(v)$ and prefix $v[1, n - \ell + 1]$. Note that in the case that $w[n - \ell, n] = w[n - \ell]^{\ell}$, the algorithm has already visited every word in $\Sigma^{*|_{P(w)}}$ with the prefix $w[1, n - \ell]$. Otherwise, as the algorithm has, by this point, visited every word of the form $w[1, n - \ell + 1] u$, for every $u \in \Sigma^{*|_{P(v[n - \ell - 1, n])}}$, it is sufficient to show that the algorithm visits every word of the form $w[1, \ell - 1] \overline{w[\ell]} u$, for every $u \in \Sigma^{* |_{P'}}, P' = P(w[n - \ell, n]) - P(\overline{w}[\ell])$.\looseness=-1
    
    Let $w'$ be the last word visited by the algorithm with the prefix $w[1,n - \ell + 1]$. Note that the first step taken by the algorithm is to determine the first position $j$ in $w'[n - \ell + 1, n]$ containing the symbol $\overline{w[\ell]}$. Therefore, by making the swap $(n - \ell, j)$, the algorithm moves to some word with a suffix in $\Sigma^{P'}$, where $P' = P(w[n - \ell, n]) - P(\overline{w[n - \ell]})$. As the algorithm can, by inductive assumption, visit every word with a suffix of length $\ell - 1$, the algorithm must also be able to visit every word with a suffix of length $\ell$, completing the proof.
\end{pf}

\begin{lem}
    \label{lem:alg_he_visits_no_node_twice}
    Let $P$ be a Parikh vector, and let $w \in \Sigma^{*|_P}$ be a word. The path output by \textsc{HamiltonEnumeration} %Algorithm \ref{alg:hamiltonian_enumeration}
    does not visit any word in $w \in \Sigma^{*|_P}$ more than once.
\end{lem}

\begin{pf}
    Note that this property holds for length 2 words. By extension, the length at most $2$ path visiting every word with the prefix $w[1, n - 2]$ does not visit the same word twice before returning to a previous call on the stack.

    Assume now that, given any input word $v \in \Sigma^{*|_P}$, the algorithm visits every word in $\Sigma^{*|_P}$ with the prefix $v[1, n - l + 1]$ without repetition, and has only visited words with this prefix. Further, assume that $P(v[n - \ell, n]) \neq (0, \ell - 1)$ or $(\ell - 1, 0)$. Then, after every such word has been visited by the path, the algorithm returns to the previous state, with the goal of enumerating every word with the prefix $v[1, n - \ell]$. As every word in $\Sigma^{*|_P}$ with the prefix $v[1, n - \ell + 1]$ has been visited, it is sufficient to show that only those words with the prefix $v[1, n - \ell] \overline{v[n - \ell + 1]}$ are enumerated. The first swap made at this state is between $\ell$ and the smallest index $j \in [n - \ell + 1, n]$ such that $v[n - \ell] \neq v[j]$, which, as the algorithm has only visited words with the prefix $v[1, n - l + 1]$, has not previously been visited. After this swap, the algorithm enumerates every word with the prefix $v[1, n - \ell - 1] \overline{v[n - \ell]}$, which, by the inductive assumption, is done without visiting the same word. Therefore, by induction. every word with the prefix $v[1, n - \ell]$ is visited by the path output by \textsc{HamiltonianEnumeration} exactly once.
\end{pf}

\begin{thm}
    \label{thm:alg_hamitonian_enumeration}
    Given a Parikh vector $P = (P_1, P_2)$ such that $P_1 + P_2 = n$, and word $w \in \Sigma^{* \vert P}$, %Algorithm \ref{alg:hamiltonian_enumeration}
    \textsc{HamiltonianEnumeration} outputs a Hamiltonian path with at most $O(\log n)$ delay between the output of each edge after $O(n \log n)$ preprocessing.
\end{thm}

\begin{pf}
    Following Lemmas \ref{lem:alg_he_visits_every_node} and \ref{lem:alg_he_visits_no_node_twice}, the path outputted by \textsc{HamiltonEnumeration} is Hamiltonian. In the preprocessing step, the algorithm constructs two balanced binary search trees $T_1$ and $T_2$. Every node in $T_1$ is labelled by some index $i_1 \in [n]$ for which $w[i_1] = 1$, and sorted by the values of the labels. Similarly, every node in $T_2$ is labelled by some index $i_2 \in [n]$ for which $w[i_2] = 2$, and sorted by the values of the labels. As each of these constructions requires at most $O(n \log n )$ time, the total complexity of the preprocessing is $O(n \log n)$.

    During each call, we have one of three cases. If either value of the Parikh vector is $0$, then the algorithm immediately returns to the last state without any output. If the Parikh vector is $(1, 1)$, then the algorithm outputs a swap between the two symbols, updates the trees $T_1$ and $T_2$, requiring at most $O(\log n)$ time, then returns to the last state. In the third case, the Parikh vector $(P_1, P_2)$ satisfies $P_1 > 0, P_2 > 0$. First, the algorithm determines the last position in the current state of the word $w$ containing the symbol $1$ and the last position containing the symbol $2$, i.e. the values $R_1 = \max_{j \in [1, n]} w[i_1] = 1$ and $R_2 = \max_{j \in [1, n]} w[i_2] = 2$. These values can be determined in $O(\log n)$ time using the trees $T_1$ and $T_2$. Using these values, the algorithm iterates through every length from $\min(T_1, T_2)$ to $n - (P_1 + P_2 - 1)$, enumerating every word in $\Sigma^{* |_{P(w)}}$ with the prefix $w[1, n - (P_1 + P_2 - 1)]$. For each $\ell \in [\min(T_1, T_2), n - (P_1 + P_2 - 1)]$, the algorithm outputs the swap $(\ell, j)$, where $j \in [n - \ell, n]$ is the largest value for which $w[j] = \overline{w[\ell]}$. After this has been output, the algorithm updates the trees $T_1$ and $T_2$. Note that both finding the value of $j$ and updating the trees require $O(\log n)$ time. After this swap, the algorithm makes the next call to \textsc{HamiltonianEnumeration}. Note that after this call, \textsc{HamiltonianEnumeration} must either return immediately to the last state or output some swap before either returning or making the next recursive call. Therefore, ignoring the time complexity of returning to a previous state in the stack, the worst case delay between outputs is $O(\log n)$, corresponding to searching and updating the trees $T_1$ and $T_2$.

    To avoid having to check each state in the stack after returning from a recursive call, the algorithm uses tail recursion. Explicitly, rather than returning to the state in the stack from which the algorithm was called, the algorithm is passed a pointer to the last state in the stack corresponding to a length $\ell$ such that some word with the prefix $w[1, n - \ell]$ has not been output. To do so, after the swap between $n - (P_1 + P_2 - 1)$ and $j$ is made, for the value $j$ as defined above, the algorithm passes the pointer it was initially given, denoted in the algorithm as $last\_state$ to the call to \textsc{HamiltonianEnumeration}, allowing the algorithm to skip over the current state during the recursion process.
\end{pf}

\noindent
\subsection{General Alphabets.}

We now show that the graph is Hamiltonian for any alphabet of size $\sigma \geq 2$. The main idea here is to build a cycle based on recursively grouping together sets of symbols. Given a Parikh vector $P = (P_1, P_2, \dots, P_{\sigma})$, our proof operates in a set of $\sigma - 1$ recursive phases, with the $i^{th}$ step corresponding to finding a Hamiltonian path in the graph $G(P_{i}, P_{i + 1}, \dots, P_{\sigma})$, then mapping this path to one in $G(P)$. The paths in $G(P_{i}, P_{i + 1}, \dots, P_{\sigma})$ are generated in turn by a recursive process. Starting with the word $w$, first, we consider the path visiting every vertex corresponding to a permutation of the symbols $i + 1, \dots, \sigma$ in $w$. Explicitly, every word $v$ in this path is of the form:\looseness=-1
$$v[i] = \begin{cases}
    w[i] & w[i] \in \{1, 2, \dots, i\}\\
    x_i  \in \{i + 1, \dots, \sigma\} & w[i] \notin \{1, 2, \dots, i\}
\end{cases},$$
where $x_i$ is some arbitrary symbol $\{i, i + 1, \dots, \sigma\}$. Further, every such word is visited exactly once.\looseness=-1

After this path is output, a single swap corresponding to the first swap in $G(P_i, (P_{i + 1} + P_{i + 2}, \dots, P_{\sigma}))$ is made, ensuring that this swap must involve some position in $w$ containing the symbol $i$. After this swap, another path visiting exactly once every word corresponding to a permutation of the symbols $i + 1, \dots, \sigma$ in $w$ can be output. By repeating this for every swap in $G(P_i, (P_{i + 1} + P_{i + 2}, \dots, P_{\sigma}))$, inserting a path visiting exactly once every word corresponding to a permutation of the symbols $i + 1, \dots, \sigma$ in $w$ between each such swap, note that every permutation of the symbols $i, i + 1, \dots, \sigma$ in $w$ is output exactly once. In other words, every word $v \in \Sigma^{*|_P}$ of the form\looseness=-1

$$v[i] = \begin{cases}
    w[i] & w[i] \in \{1, 2, \dots, i - 1\}\\
    x _i \in \{i, i + 1, \dots, \sigma\} & w[i] \in \{i, i + 1, \dots, \sigma\}
\end{cases},$$
where $x_i$ is some symbol in $\{i, i + 1, \dots, \sigma\}$. Further, each such word is visited exactly once. Using the binary alphabet as a base case, this process provides an outline of the proof of the Hamiltonicity of $G(P)$.
% A full proof of Theorem \ref{thm:general_alphabet_hamiltonian} can be found in the full version of this paper \cite{adamson2023structural}.\looseness=-1

% First, we find a Hamiltonian path performing only 2-swaps on the symbols $\sigma - 1$ and $\sigma$, ignoring all other symbols in the word.
% Once this path is found, the algorithm then makes a 2-swap using the algorithm for binary alphabets on the graph corresponding to the word formed by replacing every symbol $\sigma - 1$ and $\sigma$ with some new symbol $x_{\sigma - 1, \sigma}$.
% After making this swap, we again form a Hamiltonian path by making only 2-swaps between the symbols $\sigma - 1$ and $\sigma$.
% Note that this path now visits every word where the symbols $1, 2, \dots \sigma - 2$ are fixed as at the start of the path.

% In order to show Hamiltonicity, we introduce a set of auxiliary data structures.
% Let $x_{k}$ be a symbol corresponding to any symbol in $\Sigma$ between $\sigma - k$ and $\sigma$.
% Using this, let $w_{\ell}$ be the word formed from $w$ by removing every symbol in the set $\{1, 2, \dots, \sigma - \ell\}$, and replacing every symbol in $\{\sigma - \ell + 2, \dots, \sigma\}$ with $x_{\ell - 1}$.
% We assume that there exists a table $L_{\ell}$ mapping the indices $1,2, \dots, \vert w_{\ell} \vert$ to $1,2, \dots, n$ such that $L_{\ell}[1] < L_{\ell}[2] < \dots < L_{\ell}[\vert w_{\ell} \vert]$ and either $w_{\ell}[i] = w[L_{\ell}[i]] = \sigma - \ell$ of $w_{\ell}[i] = x_{\ell - 1}$ and $w[L_{\ell}[i]] \in \{\sigma - \ell + 2, \sigma - \ell + 3, \dots, \sigma\}$.

\begin{thm}
    \label{thm:general_alphabet_hamiltonian}
    Given an arbitrary Parikh vector $P \in \mathbb{N}^{\sigma}$, there exists a Hamiltonian path starting at every vertex $v$ in the configuration graph $G(P)$.\looseness=-1
\end{thm}

\begin{pf}
    This is formally proven using the outline above via an inductive argument with the base case of binary alphabets, the Hamiltonicity of which is proven in Theorem \ref{thm:binary_graph_is_hamiltonian}.\looseness=-1
    % for any binary Parikh vector $P$ and word $w \in \Sigma^{*|_P}$, there exists a Hamiltonian path in $G(P)$ starting at $w$.\looseness=-1

    We assume now that there exists, for any Parikh vector $p \in \mathbb{N}_0^{\ell - 1}$ and word $w \in \Sigma^{*|_P}$, there exists some Hamiltonian path in $G(P)$ starting at $w$.
    Let $q = (q_1, q_2, \dots, q_{\ell})$ be a Parikh vector, and let $v \in \Sigma^{*|_q}$ be an arbitrary word with the Parikh vector $q$.
    To construct the Hamiltonian path starting at $v$ in $G(q)$, we first form a Hamiltonian path $P_1$ in $G(q_2, q_3, \dots, q_{\ell})$ starting at the word $v'$ formed by deleting every symbol $1$ from $v$.
    We assume that we have a table $T$ such that $T[i]$ returns the index in $[1, \vert v \vert]$ for which $v'[i] = v[T[i]]$.
    To avoid any repetition, we require $T[1] < T[2] < \dots < T[\vert v' \vert]$.
    With this table, each swap $(s_1, s_2)$ in the path $P_1$ can be converted to the swap $(T[s_1], T[s_2])$ in the graph $G(q)$ swapping the same symbols in $v$ as in the reduced word $v'$.
    With this conversion, $P_1$ constructs a path in $G(q)$ visiting exactly once each word where the symbol $1$ appears only at the position $\{i \in [1, \vert v \vert] \mid v[i] = 1\}$.

    Next, we construct a Hamiltonian path $P_2$ in the graph $G(q_1, q_2 + q_3 + \dots + q_{\ell})$ starting at the word $v'$ where $v'[i] = \begin{cases}
        1 & v[i] = 1\\
        x & v[i] \neq 1
    \end{cases}$, for some new symbol $x$.
    This graph can be seen as an abstraction of $G(q)$, considering only swaps between some position labelled $1$ and any position with a different symbol.
    
    The first swap in $P_2$ is applied to the current word, however, rather than proceeding along this path, a new set of swaps is inserted corresponding to some Hamiltonian path in $G(q_2, q_3, \dots, q_{\ell})$ generated in the same manner as before.
    Again, this new path corresponds to a permutation of every symbol in the set $\{2, 3, \dots, \sigma\}$, while fixing the positions of the symbol $1$ in the word.
    This is repeated by taking a single swap from the path $P_2$, followed by a complete path corresponding to a Hamiltonian path in $G(q_2, q_3, \dots, q_{\ell})$.
    By combining these paths, the new path must visit exactly once every word in $\Sigma^{*|_{q_2, q_3, \dots, q_{\ell}}}$ where the positions of symbol $1$ are fixed, each time a word $w$ with a new permutation of the symbol $1$ is visited.
    Similarly, every word in $\Sigma^{*|_{q_1, q_2 + q_3 + \dots + q_{\sigma}}}$ is visited exactly once, corresponding to a path through every permutation of the positions of the symbols $1$ in the word.
    Therefore, the output path is Hamiltonian.
    %
    % Observe that, following Theorem \ref{thm:binary_graph_is_hamiltonian}, note that there exists a non-cyclic path covering every word in $\Sigma^{*|_P}$ where the positions of all symbols in $w$  are fixed other than $\sigma - 1$ and $\sigma$.
    % Using this observation, a Hamiltonian path on the alphabet $(\sigma - 2, \sigma - 1, \sigma)$ can be formed by constructing a Hamiltonian path in $G(p_{\sigma - 2}, p_{\sigma - 1} + p_{\sigma})$, and adding, between each swap in this path, a Hamiltonian path corresponding a Hamiltonian path in $G(p_{\sigma - 1}, p_{\sigma})$, starting at the word $w$ formed by removing all occurrences of the symbol $\sigma - 2$ from the current word.
    % As each path in $G(p_{\sigma - 1}, p_{\sigma})$ corresponds to a full permutation of every symbol $\sigma - 1, \sigma$ in the current word, and the positions of the symbols $\sigma - 2$ only change between swaps corresponding to the path in $G(p_{\sigma - 2}, p_{\sigma - 1} + p_{\sigma})$, no word is visited more than once.
    % Further, as the paths in both $G(p_{\sigma - 2}, p_{\sigma - 1} + p_{\sigma})$ and $G(p_{\sigma - 1},  p_{\sigma})$ are Hamiltonian, the path constructed for $G(p_{\sigma - 2}, p_{\sigma - 1}, p_{\sigma})$ must be Hamiltonian.
    %
    % In the general case, the same approach can be used, constructing a path $G(p)$ by combining a path in $G((p_{1}, \sum_{j \in [2, \sigma]} p_j))$ and $G((p_2, p_3, \dots, p_{\sigma}))$.
\end{pf}

\begin{thm}
    \label{thm:general_hamiltonian_algorithm}
    Given a Parikh vector $P = (P_1, P_2, \dots, P_{\sigma})$ such that \newline$\sum_{i \in [1, \sigma]} P_i = n$, there exists an algorithm outputting a Hamiltonian path in $G(P)$ with a delay of at most $O(\sigma \log n)$ between outputting each edge after $O(n \log n)$ preprocessing.\looseness=-1
\end{thm}

\begin{pf}
    See Algorithm \ref{alg:full_hamiltonain} for the full pseudocode of this algorithm.\looseness=-1

    Our algorithmic results works using the same approach as for the binary case, outlined in Theorem \ref{thm:alg_hamitonian_enumeration}. The primary difference between these algorithms is the use of a recursive proccess to enumerate the swaps between all positions containing some symbol in $(2, 3, \dots, \sigma)$, using the same algorithm to enumerate the Hamiltonian path starting with the word formed by removing every word containing the symbol $1$. In the $q^{th}$-step of the recursion, after a swap is made between some position containing the symbol $q$, and some position containing some symbol $x \in [q + 1, \sigma]$, a recursive call is made to enumerate the Hamiltonian path starting at the word formed by removing the every copy of the symbol $q$.\looseness=-1

    Our algorithm works as follows. We assume, at each step, we have the current word $w$ and $\sigma$ balanced binary search trees, $T_1, T_2, \dots, T_{\sigma}$ where $T_q$ stores all positions of the symbol $q$ in $w$. When making a call to the function, \textsc{GeneralHamiltonianEnumeration}, three parameters are passed, the Parikh vector the the current prefix being considered, $(P_1, P_2, \dots, P_{\sigma})$, the smallest symbol $q$ on which swaps are allowed, and a pointer to the return state once this proccess is exhausted. The swaps are determined in the same way as Algorithm \ref{alg:hamiltonian_enumeration}, with two key differences. Within each call, we only make swaps between positions in the word $w$ containing the symbol $q$, and positions containing any symbol in $[q + 1, \sigma]$. The swap is chosen in the same way as in the binary case, swapping the positions $R_q$, the \textbf{R}ightmost position containing $q$, and $R_p$, the \textbf{R}ightmost position containing any symbol greater than $q$. We note that, in order to determine the value of $R_p$ we determine the rightmost position of each symbol $x \in [q + 1, \sigma]$, the value of which is determined using the trees $T_{q + 1}, T_{q + 2}, \dots, T_{\sigma}$, requiring in total $\sigma \log n$ time.
    Alongside this, we store a set of positions $p_q, p_{q + 1}, \dots, p_{\sigma}$, where $p_{x} $ stores the rightmost position smaller than $\min(R_q, R_p)$ containing the symbol $x$, i.e. the value such that $p_x < \min(R_q, R_p)$ and, for every $p_x'$ where $w[p_x'] = x$, either $p_x' \geq \min(R_q, R_p)$ or $p_x' < p_x$. We use these to allow efficient update of the value $i$ denoting the largest position of the next swap.
    Once a swap has been made, two recursive calls are made. First, to \textsc{GeneralHamiltonianEnumeration} with the arguments $P(w), q + 1, current\_state()$, enumerating every permutation of the positions in $w$ containing some symbol $x \in [q + 1, \sigma]$, before returning to this call, then to \textsc{GeneralHamiltonianEnumeration}with the arguments $(P_1, P_2, \dots, P_{\sigma} - P(w[i]), q, current\_state()$, enumerating every suffix of $w$ of length $P_1 + P_2 + \dots + P_{\sigma} - 1$. Finally, once every permutation of $w[n - P_1 - P_2 - \dots - P_{\sigma} + 1, n]$ formed by swapping two positions containing distinct symbols from the set $[q, \sigma]$ has been enumerated, one last swap is made between the position $n - P_1 - P_2 - \dots - P_{\sigma}$ and the symbol at position $j$ where $j$ is the leftmost (smallest) position in $[n - (P_1 + P_2 + \dots + P_{\sigma}, n]$ that does not contain $w[n - (P_1 + P_2 + \dots + P_{\sigma}]$, i.e. the value $j$ such that $w[j] \neq w[n - (P_1 + P_2 + \dots + P_{\sigma}]$ and, $\forall j' \in [n - (P_1 + P_2 + \dots + P_{\sigma}, j - 1]$, $w[j'] = w[n - (P_1 + P_2 + \dots + P_{\sigma}]$. Once this swap has been made, the algorithm returns to the state it was passed when called, avoiding unnecessary recursion.\looseness=-1
    
    Observe that the correctness of this algorithm follows from Theorems \ref{thm:alg_hamitonian_enumeration} (laying out the correctness of our choice of swaps to make at each step) and \ref{thm:general_alphabet_hamiltonian}. Further, the worst case delay of $O(\sigma \log n)$ is due to determining the values of $R_p$ at the start of the function, and $j$ during both the main loop enumerating the permutations of the suffix, and in the final output. As each value is computed at most twice between outputs, we arrive at the stated worst case delay.\looseness=-1
\end{pf}

\begin{algorithm}[h]
    \caption{Enumeration algorithm for general alphabets}
    \label{alg:full_hamiltonain}
    \begin{algorithmic}[1]
        \State \textbf{Global Variables:}
        \State $w \in \Sigma^n$
        \State Balanced Binary Search Trees $T_1, T_2, \dots, T_{\sigma}$
        \Function{GeneralHamiltonianEnumeration}{$P \in \mathbb{N}^{\sigma}, q \in [\sigma], $, pointer $last\_state$}
            \If{$q = 2$}
                \State \textsc{HamiltonianEnumeration($(P_{\sigma - 1}, P_{\sigma})$, $CurrentState()$)}
                \State \textbf{Return To:} $last\_state$
            \ElsIf{$P_q = 0$}
                \State \textbf{Return To:} $last\_state$
            \Else
                \State $R_q \gets \max(T_q)$
                \State $R_p \gets \max_{x \in [q + 1, \sigma]} T_{x}$
                \State $i \gets \min(R_q, R_p)$
                \State $p_x \gets \max_{< i} T_x, \forall x \in [q, \sigma]$
                \While{$i > n - (P_1 + P_2 + \dots + P_{\sigma})$}
                    \State $j \gets 0$
                    \If{$w[i] = q$}
                        \State $j \gets \min_{x \in [q + 1, \sigma]} \min_{\geq i + 1} T_{x}$
                        %\State $j \gets \min_{x \in [q + 1, \sigma]} \min_{> i + 1, < n} T_{x}$
                        %\State $j \gets \min_{x \in [q + 1, \sigma]} \min_{j' \in [i + 1, n]} T_{x}$
                    \Else
                        \State $j \gets \min_{ \geq i + 1 - (\sum_{y \in [1, \sigma]} P_y)} T_{q}$
                        %\State $j \gets \min_{j' \in [i + 1 - (\sum_{y \in [1, \sigma]} P_y), n]} T_{q}$
                    \EndIf
                    \State \textbf{Output:} $(i, j)$
                    \State \textbf{Remove}$(T_{w[i], i})$, \textbf{Insert}$(T_{w[i]}, j)$
                    \State \textbf{Remove}$(T_{w[j], j})$, \textbf{Insert}$(T_{w[j]}, i)$
                    \State $w \gets w \circ (i, j)$
                    \State \textsc{GeneralHamiltonianEnumeration}$(P(w), q + 1, current\_state())$
                    \State $P' \gets (P_1, P_2, \dots, P_{q}) - P(w[i])$
                    \State \textsc{GeneralHamiltonianEnumeration}$(P', q, current\_state())$
                    \State $i, x \gets \max_{x \in [q, \sigma]} p_x$
                    \State $p_x \gets \max_{< i} T_x$
                \EndWhile
                \State $m \gets n - (P_1 + P_2 + \dots + P_{\sigma})$
                \State $j \gets \min_{x \in [q, \sigma] \setminus w[m]}\min_{\geq m} T_{x}$
                %\State $j \gets \min_{x \in [q, \sigma] \setminus w[m]}\min_{i' \in [m, n]} T_{x}$
                \State \textbf{Output} $(m, j)$
                \State \textbf{Remove}$(T_{w[j]}, j)$, \textbf{Insert}$(T_{w[j]}, m)$
                \State \textbf{Remove}$(T_{w[m]}, m)$, \textbf{Insert}$(T_{w[m]}, j)$
                \State $w \gets w \circ (m, j)$
                \State \textsc{GeneralHamiltonianEnumeration}$(P(w), q + 1, current\_state())$
                \State $P' \gets (P_1, P_2, \dots, P_n) - P(w[m])$
                \State \textsc{GeneralHamiltonianEnumeration}$(P', q, last\_state)$
            \EndIf
        \EndFunction
    \end{algorithmic}
\end{algorithm}

\section{Conclusion}
Following the work on 2-swap, the most natural step is to consider these problems for $k$-swap based on two variants with exactly $k$ and less or equal to $k$. Note that a configuration graph for exactly $k$-swap permutation might not have a single component. We also would like to point to other attractive directions of permutations on multidimensional words~\cite{CSP_undecidable} and important combinatorial objects such as necklaces and bracelets~\cite{kCentreNecklaces,rankingBracelets}. For the $2$-swap graph specifically, we leave open the problem of determining the shortest path between two given words $w$ and $v$. We conjecture that the simple greedy algorithm used to derive the upper bound in Lemma \ref{lem:diameter_upper_bound} can be used to find the shortest path between any pair of vertices.

\section{Acknowledgements}
We wish to thank Torsten M{\"u}tze for providing us with vital assistance concerning the literature of combinatorial Gray codes.
We would also like to thank the Leverhulme Trust for funding through the Leverhulme Research Centre for Functional Materials Design.
%\begin{thebibliography}{00}


%\bibitem{}

\bibliography{references}
 \bibliographystyle{elsarticle-num}
%\end{thebibliography}
\end{document}